\documentclass[12pt,twoside]{article}    
\usepackage{amssymb}
\usepackage{amsmath}
\usepackage{amsfonts}
\usepackage{latexsym}
\usepackage{graphicx}

\usepackage{url}

\newtheorem{theorem}{Theorem}[section]
\newtheorem{remark}[theorem]{Remark}
\newtheorem{lemma}[theorem]{Lemma}

\newtheorem{example}[theorem]{Example}
\numberwithin{equation}{section}
\newenvironment{proof}[1][Proof]{\textbf{#1.} }{\ \rule{0.5em}{0.5em}}

\makeatletter \@addtoreset{equation}{section} \makeatother
\begin{document}


\pagestyle{myheadings}

\markboth{\hfill {\small AbdulRahman Al-Hussein and Boulakhras Gherbal} \hfill}{\hfill
{\small Maximum principle for optimal control of FBDSDEs with jumps} \hfill }


\thispagestyle{plain}


\begin{center}
{\large \textbf{Maximum principle for optimal control of forward-backward doubly stochastic differential \\ equations with jumps}} \\
\vspace{0.7cm} {\large AbdulRahman Al-Hussein$^{a,}$\footnote{
This work is supported by the Science College Research Center at Qassim University, project no. SR-D-012-1958.},
Boulakhras Gherbal}$^{b,}$\footnote{It is also supported by the Algerian PNR project no. 8/u 07/857. }
\\
\vspace{0.2cm} {\footnotesize
{\it $^{a}$Department of Mathematics, College of Science, Qassim University, \\
 P.O.Box 6644, Buraydah 51452, Saudi Arabia \\ {\emph E-mail:} alhusseinqu@hotmail.com, hsien@qu.edu.sa \\ \smallskip
 $^{b}$Laboratory of Applied Mathematics, University of Mohamed Khider, \\ P.O.Box 145, Biskra 07000, Algeria
 \\ {\emph E-mail:} bgherbal@yahoo.fr}}
\end{center}

\begin{abstract}
In this paper we consider the maximum principle of optimal control for a stochastic
control problem. This problem is governed by a system of fully coupled multi-dimensional forward-backward doubly stochastic
differential equation with Poisson jumps. Moreover, all the coefficients appearing in this system are allowed to be random and
depend on the control variable.

We derive, in particular, sufficient conditions for optimality for this stochastic optimal control problem.
\end{abstract}



\section{Introduction}\label{sec1}
Forward-backward stochastic differential equations (FBSDEs in short) were
first studied by Antonelli in \cite{Ant93}, and since then they are
encountered in stochastic optimal control problem and mathematical finance.
For example, Xu in \cite{Xu95} studied a non-coupled continuous
forward-backward stochastic control system. Then Wu, \cite{Wu98}, studied extensively
the maximum principle for optimal control problem of fully coupled continuous
forward-backward stochastic system. We refer the reader also to \cite{Bah-Mez09}.
Peng and Wu, \cite{Pe-Wu99}, considered fully coupled continuous forward-backward
stochastic differential equations with random coefficients and applications to optimal control. A method of continuation is developed there. In this respect one can see also \cite{Yong97}.
Shi and Wu in \cite{Shi-Wu06} studied the maximum principle for fully coupled
continuous forward-backward stochastic system and provided under non-convexity assumption on the control domain necessary optimality conditions. The forward diffusion there does not contain the control variable.

Fully coupled FBSDEs with respect to Brownian motion and Poisson process were considered by Wu in \cite{W2} and
Yin and Situ in \cite{YS}.
Such equations have been shown to be very useful for example in studying linear quadratic optimal control
problems of random jumps, and also to handle nonzero-sum differential games with random jumps. The work of Wu and Wang in \cite{WW} is useful in this
respect.
In \cite{Oks-Sulem08} the authors investigated stochastic maximum
principle for non-coupled one-dimensional FBSDEs with jumps.
Meng, \cite{Ming-Max09}, considered an optimal control problem of fully
coupled forward-backward stochastic systems with Poisson jumps under partial
information. More generally, Shi in \cite{Shi2012} provided recently necessary conditions for
optimal control of fully coupled FBSDEs with random jumps.

Backward doubly stochastic differential equations were first
introduced by Pardoux and Peng in \cite{PP2}. They gave a probabilistic representation of quasi linear
stochastic partial differential equations.
In 2003 Peng and Shi, \cite{PS}, introduced fully coupled forward-backward doubly stochastic differential equations
(FBDSDEs in short). Such equations are generalizations of stochastic Hamilton systems. Existence and uniqueness of the solutions to (continuous) FBDSDEs with arbitrarily fixed time duration and under some monotone assumptions are established. Then the authors in \cite{PS} provided
also a probabilistic interpretation for the solutions of a class of quasilinear SPDEs. In this respect we refer the reader to \cite{ZS2} for an application of fully coupled FBDSDEs to provide a probabilistic formula for the solution of a quasilinear stochastic partial differential-integral equation (SPDIE in short). Another application to SPDEs can be found in \cite{Zha-Shi2012}. These are some examples to show the importance of studying FBDSDEs.

The existence and uniqueness of measurable solutions to FBDSDEs with Poisson jumps are established in \cite{ZS2} via the method of continuation. This result will be used in Section~\ref{sec2} and Section~\ref{sec3} below.
A sufficient maximum principle with partial information for a one-dimensional FBDSDE with jump with a forward equation being independent of the processes of the backward equation was studied in \cite{Jiang-Xu}. Necessary optimality conditions for FBDSDEs in \cite{Zha-Shi2012} were derived also there under non-convexity assumption on the control domain. On the other hand, in \cite{Zha-Shi2012} the authors studied the maximum principle to find necessary and sufficient conditions for optimality for a stochastic control problem governed by a continuous FBDSDE in dimension one. Within convex control domains they allow also all the coefficients of these equations to contain control variables.

The general case in particular, the maximum principle for control problems governed by a multi-dimensional discontinuous FBDSDE with its coefficients being allowed to be random and depend on the control variable and when the control domain is not convex is still an interesting incomplete research problem. In the present work we shall consider this discontinuous situation, and study, in particular, a stochastic control problem where the system is governed by a nonlinear fully coupled multi-dimensional FBDSDE with jumps as in system~(\ref{intr-syst}) below. More precisely, we shall allow both the forward and backward equations to have random jumps, and establish sufficient conditions for optimality in the form of the maximum principle. We will allow also all the coefficients appearing in our system to be random and contain control variables. Our results here are new in this respect. We will consider some relevant necessary optimality conditions for this problem in the future work.

\vspace{0.5cm}

Our system under study is the following:
\begin{eqnarray}\label{intr-syst}
\left\{
\begin{array}{ll}%
dy_{t}= b(  t,y_{t},Y_{t},z_{t},Z_{t},k_{t},v_{t})  dt+\sigma(  t,y_{t},Y_{t},z_{t},Z_{t},k_{t},v_{t})  dW_{t}\\ \hspace{2cm}
+\int_{\Theta}\varphi(  t,y_t,Y_t,z_{t},Z_{t},k_{t} ,v_{t},\rho)  \tilde{N}(  d\rho,dt)
-z_{t}\overleftarrow{dB}_{t} , \\
dY_t  = -f(  t,y_{t},Y_{t},z_{t},Z_{t},k_{t},v_{t})  dt-g(  t,y_{t},Y_{t},z_{t},Z_{t},k_{t},v_{t})
\overleftarrow{dB}_{t}\\ \hspace{2.75in}
+Z_{t}dW_{t}+\int_{\Theta}k_{t}(  \rho)  \tilde{N}(d\rho,dt) ,
\\
y_{0}=x_{0},Y_{T}=h(  y_{T})  ,
\end{array}
\right.
\end{eqnarray}
where $b,\sigma,\varphi,f,g$ and $h$ are given mappings, $(W_{t})_{t\geq0}$ and
$(B_{t})_{t\geq0}$ are independent Brownian motions taking their values respectively in $\mathbb{R}^{d}$ and $\mathbb{R}^{l},$ while $v_{\cdot}$ represents a control process and $\tilde{N}(d\rho,dt)$ is the compensated Poisson random measure associated with a Poisson point process $\eta .$ Here $T$ is a fixed positive number.

We shall be interested in minimizing the cost functional
\begin{equation}\label{eq:1.2}
J(v_{\cdot})=\mathbb{E} \big[ \int_{0}^{T} \ell(t,y_{t},Y_{t},z_{t},Z_{t},k_{t},v_{t})dt+\beta(y_{T})+\gamma(Y_{0}) \big],
\end{equation}
over the set of all admissible controls (to be described in Section~\ref{sec2} below).
%

The paper is organized as follows. In Section~\ref{sec2}, we formulate the problem and
give various assumptions used throughout the paper. In Section~\ref{sec3} we
introduce the adjoint equation of (\ref{intr-syst}), state our main theorem and give an
example to illustrate this theorem.
Section~\ref{sec4} is devoted to proving the main result.

\section{Formulation of the problem and assumptions}\label{sec2}
Let $(\Omega,\mathcal{F},\mathbb{P}) $ be a complete probability space. Let $(
W_{t}) _{t\in[ 0,T ] }$ and $( B_{t}) _{t\in[ 0,T ] }$ be two Brownian motions taking their values in $\mathbb{R}^{d}$ and $\mathbb{R}^{l}$ respectively. Let $\eta$ be a Poisson point process taking
its values in a measurable space $( \Theta,\mathcal{B}( \Theta) ) .$ We denote
by $\nu( d\rho) $ the characteristic measure of $\eta$ which is assumed to be a
$\sigma$-finite measure on $( \Theta,\mathcal{B}( \Theta) ),$ by $N (d\rho,dt)
$ the Poisson counting measure (jump measure) induced by $\eta$ with compensator $\nu( d\rho)
dt$, and by
\[
\tilde{N}( d\rho,dt) = N ( d\rho,dt) -\nu( d\rho) dt,
\]
the compensation of the jump measure $N(\cdot, \cdot)$ of $\eta .$ Hence $\nu (O) = \mathbb{E} [N(O,1)]$ for $O\in \mathcal{B}(\Theta) .$
We assume that these three processes $W , B$ and $\eta$ are mutually independent.

Let $\mathcal{N}$ denote the class of $\mathbb{P}$-null sets of $\mathcal{F}.
$ For each $t\in[ 0,T ] $, we define $\mathcal{F}_{t}:= \mathcal{F}_{t}^{W}\vee\mathcal{F}_{t,T}^{B}\vee\mathcal{F}_{t}^{\eta}$, where for any
process $\left\{  \pi_{t}\right\}  $, we set $$\mathcal{F}_{s,t}^{\pi}=\sigma(
\pi_{r}-\pi_{s};s\leq r\leq t) \vee\mathcal{N},\mathcal{F}_{t}^{\pi
}=\mathcal{F}_{0,t}^{\pi}.$$

For a Euclidean space $E,$ let $\mathcal{M}^{2}(0,T; E) $ denote the set of jointly measurable
processes $\left\{  X_{t},t\in [0 ,T] \right\}
$ taking values in $E,$  and satisfy: $X_{t}$ is $\mathcal{F}_{t}$-measurable for a.e. $t \in [0,T], $ and
\[
\mathbb{E} \big[  {\displaystyle\int_{0}^{T}} \left\vert X_{t}\right\vert_E ^{2}dt \big]
<\infty .
\]

%
Let $L_{\nu}^{2}(E)$ be the set of $\mathcal{B}( \Theta)$-measurable mapping $k$ with values in $E$
such that
\[
\left\vert \left\Vert k \right\Vert \right\vert := \big[ \int_{\Theta}\left\vert k(
\rho) \right\vert_{E}^{2}\nu( d\rho) \big]^{\frac{1}{2}} < \infty .
\]

Denote by $\mathcal{N}_{\eta}^{2}( 0,T ;E) $ to the set of processes $\{ K_t , \; t \in [0,T] \}$ that take their values in $L_{\nu}^{2}(E)$ and satisfy: $K_{t}$ is $\mathcal{F}_{t}$-measurable for a.e. $t \in [0,T], $ and
\[
\mathbb{E} \big[ {\displaystyle\int_{0}^{T}} \int_{\Theta}\left\vert K_{t}( \rho)
\right\vert_{E}^{2}\nu( d\rho) dt \big] <\infty.
\]

Finally, we set
\begin{eqnarray*}
&& \hspace{-1cm} \mathbb{M}^{2} := \mathcal{M}^{2}\left(  0,T ;\mathbb{R}^{n}\right)
\times\mathcal{M}^{2}\left( 0,T;\mathbb{R}^{m}\right)
\times\mathcal{M}^{2}\left(0,T;\mathbb{R}^{n\times l}\right)  \\ && \hspace{2.25in}
\times \mathcal{M}^{2}\left(0,T ;\mathbb{R}^{m\times d}\right)
\times\mathcal{N}_{\eta}^{2}\left(0,T ;\mathbb{R}^{m}\right)  .
\end{eqnarray*}
Then $\mathbb{M}^{2}$ is a Hilbert space
with respect to the norm $\left\Vert \cdot \right\Vert _{\mathbb{M}^{2}}$ given by
\begin{align*}
& \hspace{-0.25cm} \left\Vert \zeta_{\cdot}\right\Vert _{\mathbb{M}^{2}}^2 \\
&  \hspace{0.25cm} := \mathbb{E}\big[  \int_0^T \left\vert
y_{t}\right\vert ^{2} dt + \int_0^T \left\vert Y_{t}\right\vert ^{2}
+\int_{0}^{T}\left\Vert z_{t}\right\Vert^{2}dt+\int_{0}^{T}\left\Vert Z_{t}\right\Vert ^{2}dt+\int_{0}^{T}\left\vert \left\Vert
k_t \right\Vert \right\vert ^{2}dt \big]  ,
\end{align*}
for $\zeta_{\cdot}=\left(  y_{\cdot},Y_{\cdot},z_{\cdot},Z_{\cdot},k_{\cdot} \right)  .$

Let $U$ be a non-empty subset of $\mathbb{R}^r$. We say that
$v_{\cdot} : [0 , T]\times \Omega \rightarrow \mathbb{R}^r$ is \emph{admissible}
if $v_{\cdot} \in  \mathcal{M}^2 (0 , T; \mathbb{R}^r)$ and $v_t \in U \; \; a.e., \; \mathbb{P}-a.s.$
The set of admissible controls will be denoted by $\mathcal{U}_{ad} .$
Consider the following controlled fully coupled FBDSDE with jumps:
\begin{equation}\label{eq:2.1}
\left\{
\begin{array}{ll}%
dy_{t}= b(  t,y_{t},Y_{t},z_{t},Z_{t},k_{t},v_{t})  dt+\sigma(  t,y_{t},Y_{t},z_{t},Z_{t},k_{t},v_{t})  dW_{t}\\ \hspace{2cm}
+\int_{\Theta}\varphi(  t,y_t,Y_t,z_{t},Z_{t},k_{t} ,v_{t},\rho)  \tilde{N}(  d\rho,dt)
-z_{t}\overleftarrow{dB}_{t} , \\
dY_t  = -f(  t,y_{t},Y_{t},z_{t},Z_{t},k_{t},v_{t})  dt-g(  t,y_{t},Y_{t},z_{t},Z_{t},k_{t},v_{t})
\overleftarrow{dB}_{t}\\ \hspace{2.75in}
+Z_{t}dW_{t}+\int_{\Theta}k_{t}(  \rho)  \tilde{N}(
d\rho,dt) ,
\\
y_{0}=x_{0},Y_{T}=h(  y_{T})  ,
\end{array}
\right.
\end{equation}
where the mappings
\begin{eqnarray*}
&& b :\Omega\times[  0,T]  \times \mathbb{R}^{n}\times
\mathbb{R}^{m}\times \mathbb{R}^{n\times l}\times \mathbb{R}^{m\times d}\times L_{\nu}^{2}(  \mathbb{R}^{m})  \times\mathbb{R}^r \rightarrow\mathbb{R}^{n},\\
&& \sigma :\Omega\times[  0,T]  \times\mathbb{R}^{n}\times\mathbb{R}^{m}\times\mathbb{R}^{n\times l}\times\mathbb{R}\times L_{\nu}^{2}(  \mathbb{R}^{m})\times\mathbb{R}^r \rightarrow\mathbb{R}^{n\times d},\\
&& \varphi :\Omega\times[  0,T]  \times\mathbb{R}^{n}\times\mathbb{R}^{m}\times\mathbb{R}^{n\times l}\times\mathbb{R}^{m\times d}\times L_{\nu}^{2}(  \mathbb{R}^{m})  \times\mathbb{R}^r \times \Theta \rightarrow\mathbb{R}^{n},\\
&& f  :\Omega\times[  0,T]  \times\mathbb{R}^{n}\times\mathbb{R}^{m}\times\mathbb{R}^{n\times l}\times\mathbb{R}^{m\times d}\times L_{\nu}^{2}(  \mathbb{R}^{m})  \times\mathbb{R}^r \rightarrow
\mathbb{R}^{m},\\
&& g :\Omega\times[  0,T]  \times\mathbb{R}^{n}\times\mathbb{R}^{m}\times\mathbb{R}^{n\times l}\times\mathbb{R}\times L_{\nu}^{2}(  \mathbb{R}^{m})\times\mathbb{R}^r \rightarrow\mathbb{R}^{m\times l},\\
 && h :\Omega \times\mathbb{R}^{n}\rightarrow\mathbb{R}^{m},
\end{eqnarray*}
are measurable (further properties to be introduced later in this section) and $v_{\cdot} \in \mathcal{U}_{ad} .$ Given a full-rank $m\times n$ matrix $R$  of real indices, we assume that $h$ is defined, for $(\omega , x) \in \Omega \times\mathbb{R}^{n} ,$ by $h(\omega , x) := c \,  R x + \xi (\omega) , $ where $c \neq 0$ is a constant and $\xi$ is a fixed arbitrary element of $L^{2}(  \Omega,\mathcal{F}_{T},\mathbb{P};\mathbb{R}^{m}) .$

Note that the integral with respect to $\overleftarrow{dB}$ is a ``backward" It\^{o} integral,
while the integral with respect to $dW$ is a standard ``forward"
It\^{o} integral. We refer the reader to \cite{Nua-Pard} for more details on such integrals, which are particular cases of the It\^{o}-Skorohod stochastic integral.

A solution of (\ref{eq:2.1}) is a quintuple $( y,Y,z,Z,k )$ of stochastic processes
such that $( y,Y,z,Z,k )$ belongs to $\mathbb{M}^{2}$ and satisfies the following FBDSDE:
\begin{eqnarray*}
\left\{
\begin{array}{ll}%
\!\! y_{t}= x_{0} + \int_0^t b(  s,y_{s},Y_{s},z_{s},Z_{s},k_{s}
,v_{s})  ds  + \int_0^t \sigma(  s,y_{s},Y_{s},z_{s},Z_{s},k_{s},v_{s})  dW_{s} \\ \hspace{2cm}
+\int_0^t \int_{\Theta}\varphi(  s,y_{s},Y_{s},z_{s},Z_{s},k_{s},v_{s},\rho)  \tilde{N}(  d\rho,ds)
- \int_0^t z_{s}\overleftarrow{dB}_{s} , \\ \\
\!\! Y( t)  = h(y_T) + \int_t^T f(  s,y_{s},Y_{s},z_{s},Z_{s},k_{s},v_{s})  ds \\ \hspace{2in} + \int_t^T g(  s,y_{s},Y_{s},z_{s},Z_{s},k_{s},v_{s})
\overleftarrow{dB}_{s} \\ \hspace{2in}
- \int_t^T Z_{s}dW_{s}- \int_t^T \int_{\Theta}k_{s}(  \rho)  \tilde{N}(
d\rho,ds) , \; \; t \in [0,T] .
\end{array}
\right.
\end{eqnarray*}

Define the cost functional by:
\begin{equation}\label{eq:2.2}
J(v_{\cdot}):=\mathbb{E} \big[ \int_{0}^{T} \ell(t,y_{t},Y_{t},z_{t},Z_{t},k_{t},v_{t})dt+\beta(y_{T})+\gamma(Y_{0}) \big],\;v_{\cdot} \in\mathcal{U}_{ad} ,
\end{equation}
where
\begin{align*}
\beta &  :\mathbb{R}^{n}\rightarrow\mathbb{R},\\
\gamma &  :\mathbb{R}^{m}\rightarrow\mathbb{R},\\
\ell &  :\Omega\times\lbrack0,T]\times\mathbb{R}^{n}\times\mathbb{R}^{m}%
\times\mathbb{R}^{n\times l}\times\mathbb{R}^{m\times d}\times L_{\nu}^{2}(\mathbb{R}^{m})\times\mathbb{R}^{r}\rightarrow\mathbb{R}
\end{align*}
are measurable functions such that (\ref{eq:2.2}) is defined.

Now the control problem of system~(\ref{eq:2.1}) is to minimize $J$ over $\mathcal{U}_{ad} .$ In this case we
say that $u_{\cdot}\in\mathcal{U}_{ad}$ is an \emph{optimal control} if
\begin{equation}\label{eq:2.3}
J( u_{\cdot} ) =\inf_{v_{\cdot} \in\mathcal{U}_{ad}} J( v_{\cdot}) .
\end{equation}

Let us set the following notations:
\begin{eqnarray*}
\zeta &=& (  y,Y,z,Z,k) \in\mathbb{R}^{n+m+n\times
l+m\times d}\times L_{\nu}^{2}(  \mathbb{R}^{m})  , \\
A(t,\zeta , v) &=& (  - R^{\ast}f,Rb, - R^{\ast}g, R\sigma, R\varphi) (t,\zeta , v)  ,
\\
\left\langle A,\zeta\right\rangle  &=& - \left\langle y,R^{\ast}f\right\rangle
+\left\langle Y,Rb\right\rangle - \left\langle z,R^{\ast}g\right\rangle
+\left\langle Z,R\sigma\right\rangle +\left\langle \left\langle k ,R\varphi
\right\rangle \right\rangle ,
\end{eqnarray*}
where%
\begin{eqnarray*}
R^{\ast}g  &=& (  R^{\ast}g_{1}\cdots R^{\ast}g_{l})  ,R\sigma=(
R\sigma_{1}\cdots R\sigma_{d})  , \ldots ,
\end{eqnarray*}
by using the columns $\{g_1, \ldots, g_l\}$ and $\{\sigma_1, \dots, \sigma_d\}$ of $g$ and $\sigma$ respecively, and
\begin{eqnarray*}
\left\langle \left\langle k ,R\varphi\right\rangle \right\rangle (t,\zeta , v)
&=& \int_{\Theta}\left\langle k(  \rho)  ,R\varphi(
t,\zeta,v,\rho)  \right\rangle \nu(  d\rho)  .
\end{eqnarray*}

The following assumptions will be our main assumptions in the paper. We
shall mimic similar assumptions from the literature (e.g. \cite{Zha-Shi2012}) for this purpose.
\begin{itemize}
\item{(A1)} \, ${ \forall \, \zeta=(  y,Y,z,Z,k)
,\bar{\zeta}=(  \bar{y},\bar{Y},\bar{z},\bar{Z},\bar{k})  \in
\mathbb{R}^{n+m+n\times l+m\times
d}\times L_{\nu}^{2}(  \mathbb{R}^{m}) , } \; \\ \forall \, t\in [   0,T ]   , \; \forall \, v \in \mathbb{R}^r ,$
\begin{eqnarray*}
&& \hspace{-1cm} \left
\langle A(t,\zeta , v)  - A(  t,\bar{\zeta} , v)
,\zeta-\bar{\zeta}\right\rangle
\leq -\mu_1 (\left\vert R(  y-\bar{y})  \right\vert^{2} + \left\vert
R^{\ast} (  Y-\bar{Y})  \right\vert ^{2} )\\
&& \hspace{2cm} - \mu_2 ( \| R (  z-\bar{z})  \|^{2}+\left\Vert
R^{\ast}(  Z-\bar{Z})  \right\Vert ^{2}+\left\vert \left\Vert
R^{\ast}(  k -\bar{k} )  \right\Vert \right\vert^{2}),
\end{eqnarray*}
and%
\[
 c>0,
\]
or
\item{(A1)'}
\begin{eqnarray*}
&& \hspace{-1cm} \left\langle A(  t,\zeta , v)  - A(  t,\bar{\zeta} , v)
,\zeta-\bar{\zeta}\right\rangle \geq \mu_1 (\left\vert R(  y-\bar{y})  \right\vert^{2}+\left\vert
R^{\ast} ( Y-\bar{Y})  \right\vert ^{2} ) \\
&& \hspace{2cm} + \mu_2 ( \| R (  z-\bar{z})  \|^{2}+\left\Vert
R^{\ast}(  Z-\bar{Z})  \right\Vert ^{2}+\left\vert \left\Vert
R^{\ast}(  k - \bar{k})  \right\Vert \right\vert^{2}),
\end{eqnarray*}
and%
\begin{eqnarray*}
c<0,
\end{eqnarray*}
where $\mu_1$ and  $\mu_2$ are nonnegative constants with $\mu_1 + \mu_2 > 0 .$ Moreover $\mu_1 > 0$ (resp. $\mu_2 > 0$) when $m > n$ (resp. $n > m$).
\item{(A2)} For each $\zeta\in\mathbb{R}^{n+m+n\times
l+m\times d}\times L_{\nu}^{2}(  \mathbb{R}^{m})  , v \in \mathbb{R}^r , \, A(t,\zeta , v)  \in \mathbb{M}^{2} .$
\item{(A3)} We assume that
\begin{eqnarray*}
\hspace{-1.1cm} \left\{
\begin{array}{ll}
\hspace{-0.15cm} (i) & \hspace{-0.5cm} \text{the mappings} \;  f,b,g,\sigma,\varphi, \ell \; \text{ are
continuously differentiable with respect to } \;   \\ & \hspace{-0.5cm} (y,Y,z,Z,k,v) , \; \text{and} \; \beta \; \text{ and}  \; \gamma \; \text{are
continuously differentiable with respect}   \\ & \hspace{-0.25cm} \text{to} \; y \; \text{and} \; Y , \; \text{respectively}, \\
\hspace{-0.15cm} (ii) & \hspace{-0.25cm} \text{the derivatives of }  f,b,g,\sigma,\varphi \; \text{with respect to the above arguments are} \\ & \hspace{-0.25cm} \text{bounded}, \\
\hspace{-0.20cm} (iii) & \hspace{-0.3cm} \text{the derivatives of } \ell \text{\ are bounded by }  C (1+\left\vert y\right\vert +\left\vert Y\right\vert +\| z\|
+\| Z\| +\left\vert \left\Vert k\right\Vert \right\vert
)  , \\
\hspace{-0.15cm} (iv) & \hspace{-0.255cm} \text{the derivatives of } \; \beta \; \text{and} \; \gamma \; \text{are bounded by} \; C \, (  1+\left\vert y\right\vert )  \; \text{ and } \; C \, (1+\left\vert Y\right\vert ) \\ & \hspace{-0.25cm} \text{respectively},
\end{array}
\right.
\end{eqnarray*}
for some constant $C > 0.$
\end{itemize}
\begin{remark}
The condition $c > 0$ in (A1) guarantees the following monotonicity condition of the mapping $h$:
\[
\left\langle h(  y)  -h(  \bar{y})  ,R(  y-\bar
{y})  \right\rangle \geq c \, | R(  y-\bar{y}) |^2 , \; \; \forall \; y , \bar{y} \in \mathbb{R}^n.
\]
The same thing happens also for $c <0$ in (A1)'.
\end{remark}

\bigskip

The following theorem is concerned with the existence and uniqueness of the solution of (\ref{eq:2.1}).
\begin{theorem}\label{Propn 2.1} For any given admissible control
$v_{\cdot} ,$ if assumptions (A1)--(A3) (or (A1)', (A2), (A3)) hold, then (\ref{eq:2.1}) has a unique solution.
\end{theorem}

Our assumptions in this theorem satisfy the assumptions of the corresponding result in \cite{ZS2}, so the proof of this theorem can be gleaned from \cite{ZS2}.

\section{Adjoint equations and the maximum principle}\label{sec3}
Suppose that (A1)--(A3) hold. We want to introduce the adjoint equations of FBDSDE~(\ref{eq:2.1}) and then present our main result of the maximum principle for our optimal control problem governed by the FBDSDE with jumps (\ref{eq:2.1}). To this end, let us begin by defining the Hamiltonian $H$ from $[0,T]\times \Omega \times\mathbb{R}^{n} \times\mathbb{R}^{m}\times\mathbb{R}^{n\times l}\times\mathbb{R}^{m\times d}\times
L_{\nu}^{2}(\mathbb{R}^{m})\times\mathbb{R}^{r}\times
\mathbb{R}^{m}\times\mathbb{R}^{n}\times\mathbb{R}^{m\times l}\times\mathbb{R}^{n\times d}\times
L_{\nu}^{2}(\mathbb{R}^{n})$ to $\mathbb{R}$ by the formula:
\begin{eqnarray}\label{eq:3.1}
&& \hspace{-2cm} H(  t,y,Y,z,Z,k,v,p,P,q,Q,V) := \left\langle p,f(  t,y,Y,z,Z,k ,v)
\right\rangle \nonumber \\
&& - \left\langle P,b(  t,y,Y,z,Z,k,v)  \right\rangle +\left\langle q,g(  t,y,Y,z,Z,k,v)  \right\rangle \nonumber \\
&& -\left\langle
Q,\sigma(  t,y,Y,z,Z,k,v)  \right\rangle - \ell(  t,y,Y,z,Z,k ,v) \nonumber \\
&& - \int_{\Theta} \left\langle V(  \rho), \varphi(  t,y,Y,z,Z,k,v,\rho) \, \right\rangle \nu(  d\rho) .
\end{eqnarray}

Let $v_{\cdot}$ be an arbitrary element of $\mathcal{U}_{ad}$ and
$\{(y_t,Y_t,z_t,Z_t,k_t) , \; t \in [0, T] \}$ be the corresponding solution of (\ref{eq:2.1}). The adjoint equations of our FBDSDE with jumps~(\ref{eq:2.1}) are
\begin{eqnarray}\label{eq:3.2}
\left\{
\begin{array}{ll}
dp_{t}=H_{Y}dt+H_{Z}dW_{t}-q_{t}\overleftarrow{dB}_{t}+\int_{\Theta}H_{k}\tilde{N}(
d\rho,dt) , \\
dP_{t}=H_{y}dt+H_{z}\overleftarrow{dB}_{t}+Q_{t}dW_{t}+\int_{\Theta}V_{t}(  \rho)
\tilde{N}(  d\rho,dt) , \\
p_{0}=-\gamma_{Y}(  Y_{0})  ,P_{T}=-c \, R^{\ast}  \,
p_{T}+\beta_{y}(  y_{T})  ,
\end{array}
\right.
\end{eqnarray}
where $H_y$ is the gradient $\nabla_y H ( t , y , Y_t, z_t, Z_t, k_t , v_t,p_t,P_t,q_t,Q_t,V_t) \in \mathbb{R}^n , \ldots $ etc.
Let us say some thing more about this system (\ref{eq:3.2}).
\begin{theorem}\label{thm:existence for adjoint eqn}
Under (A1)--(A3) there exists a unique solution $(p,P,q,Q,V)$ of the adjoint equations~(\ref{eq:3.2}) (in $\tilde{\mathbb{M}}^2 : =\mathcal{M}^{2}\left(  0,T ;\mathbb{R}^{m}\right)\times\mathcal{M}^{2}\left( 0,T;\mathbb{R}^{n}\right)
\times\mathcal{M}^{2}\left(0,T;\mathbb{R}^{m\times l}\right)
\times \mathcal{M}^{2}\left(0,T ;\mathbb{R}^{n\times d}\right)
\times\mathcal{N}_{\eta}^{2}\left(0,T ;\mathbb{R}^{n}\right)$) .
\end{theorem}
\begin{proof}
This system (\ref{eq:3.2}) can be rewritten as in the following system:
\[
\left\{
\begin{array}{ll}
dp_{t}= (  f_{Y}^{\ast}p_{t}-b_{Y}^{\ast}P_{t}+g_{Y}^{\ast}q_{t}-\sigma_{Y}^{\ast}Q_{t}%
-\int_{\Theta}\varphi_{Y}^{\ast}V_{t}(\rho)\nu(  d\rho)
-\ell_{Y})  dt \\  \hspace{1cm}
+(  f_{Z}^{\ast}p_{t}-b_{Z}^{\ast}P_{t}+g_{Z}^{\ast}q_{t}-\sigma_{Z}^{\ast}Q_{t}-
\int_{\Theta} \varphi_{Z}^{\ast}V_{t}(\rho)\nu(  d\rho)  -\ell_{Z})
dW_{t}-q_{t}\overleftarrow{dB}_{t} \\  \hspace{1.25cm}
+\int_{\Theta}(  f_{k}^{\ast}p_{t}-b_{k}^{\ast}P_{t}+g_{k}^{\ast}q_{t}-\sigma_{k}^{\ast}Q_{t}-
(\varphi_{k}^{\ast} V_{t})(\rho) - \ell_{k})
\tilde{N}(  d\rho,dt) , \\ \\
dP_{t}= (  f_{y}^{\ast}p_{t}-b_{y}^{\ast}P_{t}+g_{y}^{\ast}q_{t}-\sigma_{y}^{\ast}Q_{t}%
-\int_{\Theta}\varphi_{y}^{\ast}V_{t}(\rho)\nu(  d\rho)
-\ell_{y})  dt \\  \hspace{1.25cm}
+(  f_{z}^{\ast}p_{t}-b_{z}^{\ast}P_{t}+g_{z}^{\ast}q_{t}-\sigma_{z}^{\ast}Q_{t}-
\int_{\Theta}\varphi_{z}^{\ast}V_{t}(\rho)\nu(  d\rho)  - \ell_{z})
\overleftarrow{dB}_{t}+Q_{t}dW_{t}\\  \hspace{3.75in}
+\int_{\Theta}V_{t}(  \rho)  \tilde{N}(  d\rho,dt) , \\ \\
p_{0} = - \gamma_{Y}(  Y_{0})  , P_{T} = -c \, R^{\ast}  \,
p_{T}+\beta_{y}(  y_{T}) ,
\end{array}
\right.
\]
which is a linear FBDSDE with jumps. Here $f_{y}^{\ast} \in \mathbb{R}^{n \times m} $ is the adjoint of the Fr\'echet derivative $D_y f ( t , y , Y_t, z_t, Z_t, k_t , v_t) \in \mathbb{R}^{m\times n}$ of $f ( t , \cdot , Y_t, z_t, Z_t, k_t , v_t): \mathbb{R}^n \rightarrow \mathbb{R}^m$ at $y , \ldots$ etc., and $\beta_{y}(  y_{T})$ is the gradient $\nabla_y \beta ( y_{T}) \in \mathbb{R}^n .$

Thanks to assumptions (A1)--(A3) this latter linear FBDSDE satisfy easily (A1)', (A2) and (A3). Thus the desired result follows from Theorem~\ref{Propn 2.1}.
\end{proof}

\bigskip

Now our main theorem is the following.
\begin{theorem}\label{Theorem.3.1.}
Assume that (A1)--(A3) hold. Given $u_{\cdot} \in \mathcal{U}_{ad}$, let $(y,Y,z,Z,k)$ and $(p,P,q,Q,V)$ be the corresponding solutions of the FBDSDEs (\ref{eq:2.1}) and (\ref{eq:3.2}) respectively. Suppose that the following assumptions hold: \\
(i) \, $\beta$ and $\gamma$ are convex, \\
(ii) for all $t \in [0 , T] , \; \mathbb{P}$\,-\,a.s., the function $H( t, \cdot, \cdot, \cdot, \cdot, \cdot, \cdot,p_t,P_t,q_t,Q_t,V_t) $ is concave, \\
(iii) we have
\begin{eqnarray}\label{eq:3.3}
&& \hspace{-1.5cm} H(  t,y_{t},Y_{t},z_{t},Z_{t},k_{t},u_{t},p_{t},P_t, q_{t},Q_{t},V_{t}) \nonumber \\
&& = \max_{v \in U} H(  t,y_{t},Y_{t},z_{t},Z_{t},k_{t},v,p_{t},P_t, q_{t},Q_{t},V_{t})  ,
\end{eqnarray}
for a.e, $\mathbb{P}$\,-\,a.s.

Then $(y,Y,z,Z,k,u_{\cdot})$ is an optimal solution of the
control problem (\ref{eq:2.1})--(\ref{eq:2.3}).
\end{theorem}
\begin{remark}
Condition (A1) assumed in this theorem and the lemmas that follow is only needed to guarantee the existence and uniqueness of the solutions of (\ref{eq:2.1}) and (\ref{eq:3.2}), and so if we can get these unique solutions without assuming (A1) there will not any necessity to assume (A1) in this theorem.
\end{remark}

The proof of this theorem will be established in Section~\ref{sec4}. Now to illustrate this theorem let us present an example.
\begin{example}\label{ex1}
Let $(\Theta , \mathcal{B}(\Theta)) = ([0,1], \mathcal{B}([0,1])) .$ Let $\tilde{N}(d\rho,dt)$ be a compensated Poisson random measure, where $(t,\rho)\in[0,1]\times[0,1]$. Recall that $\mathbb{E}[\tilde{N}(d\rho,dt)^{2}]=\nu(  d\rho)dt$ is required to be a finite Borel measure such that $\int_{[0,1]}\rho^{2} \nu( d\rho)<\infty$.
Let the controls domain be $U=[-1,1]$.
Consider the following stochastic control system:
\begin{eqnarray}\label{eq:3.4}
\left\{
\begin{array}{ll}
dy_{t}=(1+t)v_{t}dt+(-z_{t}+Z_{t}+ \int_{\Theta} k_{t}(\rho) \nu(  d\rho) +v_{t})dW_{t}-z_{t}\overleftarrow{dB}_{t} \\
\hspace{9cm} - \, v_{t}\int_{[0,1]}\rho\tilde{N}(
d\rho,dt) , \\ \\
dY_{t}=-(t-4)v_{t}dt-\frac{3}{2}(z_{t}+Z_{t}+ \int_{\Theta} k_{t}(\rho) \nu(  d\rho) +v_{t})\overleftarrow{dB}_{t}+Z_{t}dW_{t}\\
\hspace{9cm} + \, \int_{\Theta}k_{t}(  \rho)\tilde{N}(  d\rho,dt) , \\
y_{0}=Y_{1}=x \in \mathbb{R} , \; t \in (0,1)   ,
\end{array}
\right.
\end{eqnarray}
where $W, B$ are Brownian motions in $\mathbb{R} ,$ and $W, B$ and $\tilde{N}$ are mutually independent.
Consider also a cost functional given for $v_{\cdot} \in\mathcal{U}_{ad}$ by
\begin{eqnarray}\label{eq:3.5}
J(v_{\cdot})=\frac{1}{2}\; \mathbb{E}\big[ \int_{0}^{1}(y_{t}^{2}+Y_{t}^{2}+z_{t}^{2}+Z_{t}^{2}+ \int_{\Theta} k_{t}^2 (\rho) \nu(  d\rho)
+v_{t}^{2})dt + y_{1}^{2} + Y_{0}^{2} \big] .
\end{eqnarray}
We define the value function by
\begin{equation}
J( u^{*}_{\cdot} ) =\inf_{v_{\cdot} \in\mathcal{U}_{ad}} J( v_{\cdot}) . \label{eq:3.6}%
\end{equation}

This system (\ref{eq:3.4}) can be related to the one in (\ref{eq:2.1}) by setting the following mappings:
\begin{eqnarray*}
&&  b(  t,y_{t},Y_{t},z_{t},Z_{t},k_{t}  ,v_{t}
)=(1+t)v_{t},\\
&&  \sigma(  t,y_{t},Y_{t},z_{t},Z_{t}
,k_{t},v_{t} )= -z_{t}+Z_{t} +\int_{\Theta} k_t (\rho) \nu(  d\rho) +v_{t} ,\\
&& \varphi(  t,y_t,Y_t,z_{t},Z_{t},k_{t},v_{t},\rho)=-v_{t}\rho,\\
&& f(  t,y_{t},Y_{t},z_{t},Z_{t},k_{t}
,v_{t})=(t-4)v_{t},\\
&& g(  t,y_{t},Y_{t},z_{t},Z_{t},k_{t}
,v_{t})=\frac{3}{2}(z_{t}+Z_{t}+ \int_{\Theta} k_t (\rho) \nu(  d\rho) +v_{t}), \\
&& h( y_1)=y_1,  \, \text{i.e.}  \; c=1, \xi=0 , \\
&& \ell(  t,y_{t},Y_{t},z_{t},Z_{t},k_{t},v_{t})=\frac{1}{2}(y_{t}^{2}+Y_{t}^{2}+z_{t}^{2}+Z_{t}^{2}+ \int_{\Theta} k_{t}^2 (\rho) \nu(  d\rho)+v_{t}^{2}) , \\
&& \beta(y_t)=\gamma(y_t)=\frac{1}{2}y_t^{2} .
\end{eqnarray*}


Letting $u_{\cdot}\equiv0 ,$ we find from the construction of FBDSDEs with jumps (as for instance in \cite{ZS2}) that the corresponding solution  $(y_{t},Y_{t},z_{t},Z_{t},k_{t})$ of (\ref{eq:3.4}) equals $(x,x,0,0,0), $ for all $t \in [0,1] .$

Next notice that the adjoint equations of (\ref{eq:3.4}) are
\begin{eqnarray}\label{eq:3.7}
\left\{
\begin{array}{ll}
dp_{t}=-Y_{t}dt+(\frac{3}{2}q_{t}-Q_{t}-Z_{t})dW_{t}-q_{t}\overleftarrow{dB}_{t} \\ \hspace{2in} +\int_{\Theta}(\frac{3}{2}q_{t}-Q_{t} -k_{t}(\rho) )\tilde{N}(d\rho,dt) , \\
dP_{t}=-y_{t}dt+(\frac{3}{2}q_{t}+Q_{t}-z_{t})\overleftarrow{dB}_{t}+Q_{t}dW_{t}+\int_{\Theta}V_{t}( \rho)\tilde{N}(d\rho,dt) , \\
p_{0}=-x  ,P_{1}= - p_{1}+ x , \; t \in (0,1) .
\end{array}
\right.
\end{eqnarray}
Since $p_{0}$ is deterministic, then so is $p_{t} .$ Hence
\begin{eqnarray*}
p_{t}=p_{0}-\int_{0}^{t}Y_{t}dt=p_{0}-x\int_{0}^{t}dt=-x-xt=-x(1+t) .\\
\end{eqnarray*}
Thus $P_{1}$ is deterministic since:
\begin{eqnarray*}
P_{1} = - p_{1}+ x = 3x .
\end{eqnarray*}
It follows similarly that
\begin{equation*}
P_{t}=P_{1}+\int_{t}^{1}y_{t}dt=3x+x(1-t)=x(4-t) .
\end{equation*}

In particular, $(p_{t},P_{t},q_{t},Q_{t},V_{t})\equiv(-x(1+t),x(4-t),0,0,0)$ is the unique solution of (\ref{eq:3.7}).
These facts show that the Hamiltonian attains an explicit formula:
\begin{eqnarray*}
&& \hspace{-1cm} H(  t,y_{t},Y_{t},z_{t},Z_{t},k_{t}  ,v,p_{t}
,P_t,q_{t},Q_{t},V_{t}) =p_{t}(t-4)v-P_{t}(1+t)v \\
&& + \, \frac{3}{2}q_{t}(z_{t}+Z_{t}+\int_{\Theta} k_t (\rho) \nu(  d\rho) +v  )-Q_{t}(-z_{t}+Z_{t}+\int_{\Theta} k_t (\rho) \nu(  d\rho) +v  )\\
&& -\int_{\Theta} v \, \rho \, V_{t}(\rho) \, \nu(d \rho) -\frac{1}{2}(y_{t}^{2}+Y_{t}^{2}+z_{t}^{2}+Z_{t}^{2}+\int_{\Theta} k_t^2 (\rho) \nu(  d\rho) +v^{2}) \\
&& =-x(1+t)(t-4)v-x(1+t)(4-t)v-\frac{1}{2}v^{2}-\frac{1}{2}x^{2}-\frac{1}{2}x^{2}\\
&& =-\frac{1}{2}v^{2}-x^{2} , \; v \in U .
\end{eqnarray*}
Hence
\begin{eqnarray*}
&& \hspace{-1cm} H(  t,y_{t},Y_{t},z_{t},Z_{t},k_{t}  ,v,p_{t}%
,P_t,q_{t},Q_{t},V_{t})-H(  t,y_{t},Y_{t},z_{t},Z_{t},k_{t}
,u_{t},p_{t},P_t,q_{t},Q_{t},V_{t}) \\ && \hspace{1.5cm} = \, - \, \frac{1}{2}v^{2}-x^{2}+\frac{1}{2} u_{t}^{2}+x^{2} = -\frac{1}{2}v^{2}
\leq 0, \; \forall v\in U, \; a.e \; t, \;  \mathbb{P}-a.s.
\end{eqnarray*}
As a result, condition (iii) of Theorem~\ref{Theorem.3.1.} holds here for $u_{\cdot} = 0 .$ Furthermore, all other conditions of Theorem~\ref{Theorem.3.1.} can be verified easily. Consequently, $$(y,Y,z,Z,k,u_{\cdot})\equiv(x,x,0,0,0,0)$$ is an optimal solution of the control problem (\ref{eq:3.4})--(\ref{eq:3.6}).
\end{example}

For more applications of the theory of fully coupled FBDSDEs particularly in providing a
probabilistic formula for the solution of a quasilinear SPDIE we refer the reader to \cite[P. 15]{ZS2}.

\section{Proofs}\label{sec4}
In this section we shall establish the proof of Theorem~\ref{Theorem.3.1.}.
Let us recall first the following lemma.
\begin{lemma}[Integration by parts]\label{Lemma 4.1.}
Let $(\alpha,\widehat{\alpha})\in
\big[ \mathcal{S}^{2}(0,T;\mathbb{R}^{n}) \big]^{2},(\beta,\widehat{\beta})\in\big[ \mathcal{M}^{2}(0,T;\mathbb{R}^{n}) \big]^{2},
(\gamma,\widehat{\gamma})\in\big[ \mathcal{M}^{2}(0,T;\mathbb{R}^{n\times k}) \big]^{2},(\delta,\widehat{\delta})
\in \big[ \mathcal{M}^{2}(0,T;\mathbb{R}^{n\times d}) \big]^{2}$, and $(K,\widehat{K})\in
 \big[\mathcal{N}_{\eta}^{2}(0,T;\mathbb{R}^{m})\big]^{2} .$ Assume that
\[
\alpha_{t}=\alpha_{0}+\int_{0}^{t}\beta_{s}ds+\int_{0}^{t}\gamma_{s}%
\overleftarrow{dB}_{s}+\int_{0}^{t}\delta_{s}dW_{s}+\int_{0}^{t}\int_{\Theta}K_{s}(\rho)\tilde{N}(d\rho,ds),
\]
and
\[
\widehat{\alpha}_{t}=\widehat{\alpha}_{0}+\int_{0}^{t}\widehat{\beta}_{s}ds+\int_{0}^{t}\widehat{\gamma}_{s}
\overleftarrow{dB}_{s}+\int_{0}^{t}\widehat{\delta}_{s}dW_{s}+\int_{0}^{t}\int_{\Theta}\widehat{K}_{s}(\rho)\tilde{N}(d\rho,ds),
\]
for $t \in [0,T] .$
Then
\begin{eqnarray*}
&& \langle \alpha_{T},\widehat{\alpha}_{T}\rangle =\langle \alpha_{0},\widehat{\alpha}_{0}\rangle+\int_{0}^{T}\left\langle \alpha_{t},d\widehat{\alpha}_{t}\right\rangle+\int
_{0}^{T}\left\langle \widehat{\alpha}_{t} , d\alpha_{t} \right\rangle+\int_{0}^{T} d\left\langle \alpha,\widehat{\alpha}\right\rangle_{t} .
\end{eqnarray*}%
\begin{eqnarray*}
&& \hspace{-0.5cm} \mathbb{E}\big[
\langle \alpha_{T},\widehat{\alpha}_{T}\rangle \big] =  \mathbb{E}\big[ \langle \alpha_{0},\widehat{\alpha}_{0}\rangle \big] + \mathbb{E} \big[ \int_{0}^{T}\left\langle \alpha_{t},d\widehat{\alpha}_{t}\right\rangle \big] + \mathbb{E} \big[ \int_{0}^{T}\left\langle
\widehat{\alpha}_{t} , d\alpha_{t} \right\rangle \big]
\\ && \hspace{0.5cm} -\mathbb{E}\big[ \int_{0}^{T}\langle\gamma_{t},\widehat{\gamma}_{t}\rangle dt \big]+\mathbb{E}\big[ \int_{0}^{T}\langle\delta_{t},\widehat{\delta}_{t}\rangle dt \big]+\mathbb{E}\big[ \int_{0}^{T}\int_{\Theta}\langle K_{t}(\rho) ,\widehat{K}_{t}(\rho) \rangle \nu(  d\rho) dt \big].
\end{eqnarray*}
\end{lemma}

This lemma can be deduced directly from It\^{o}'s formula with jumps (see e.g. \cite{SL} and \cite{Situ}).

\bigskip

We now prove Theorem~\ref{Theorem.3.1.}. We start with two lemmas.

\begin{lemma}\label{proof-lemma1}
Assume (A1)--(A3). Let $v_{\cdot}$ be an arbitrary element of $\mathcal{U}_{ad},$ and let $(y^{v_{\cdot}} ,Y^{v_{\cdot}} , z^{v_{\cdot}} , Z^{v_{\cdot}} , k^{v_{\cdot}})$ be the corresponding solution of (\ref{eq:2.1}). Then we have
\begin{eqnarray}\label{eq:4.1}
&& \hspace{-0.75cm} J(  v_{\cdot})  - J(  u_{\cdot}) \geq\mathbb{E} \big[ \left\langle   P_{T},
y_{T}^{v_{\cdot}}-y_{T}\right\rangle \big]   +\mathbb{E} \big[ \left\langle  c \, R^{\ast}
p_{T},  y_{T}^{v_{\cdot}}-y_{T}\right\rangle  \big] -\mathbb{E} \big[ \left\langle  p_{0},
Y_{0}^{v_{\cdot}}-Y_{0}\right\rangle   \big] \nonumber \\
&& \hspace{0.5cm} + \; \mathbb{E} \big[ \int_{0}^{T} \big(   \ell(  t,y_{t}^{v_{\cdot}},Y_{t}^{v_{\cdot}},z_{t}^{v_{\cdot}},
Z_{t}^{v_{\cdot}},k_{t}^{v_{\cdot}}  ,v_{t})  -\ell(t,y_{t},Y_{t},z_{t},Z_{t},k_{t}  ,u_{t})   \big)
dt \big] .
\end{eqnarray}
\end{lemma}
\begin{proof}
From (\ref{eq:2.2}) we get
\begin{eqnarray*}
&& J(  v_{\cdot})  - J(  u_{\cdot}) = \mathbb{E} \big[  \beta(y_{T}^{v_{\cdot}})  -\beta(  y_{T}) \big]  +\mathbb{E} \big[
\gamma(  Y_{0}^{v_{\cdot}})  -\gamma(  Y_{0})  \big] \\
&& \hspace{1.5cm} + \; \mathbb{E} \big[ \int_{0}^{T} \big( \,  \ell(  t,y_{t}^{v_{\cdot}},Y_{t}^{v_{\cdot}},z_{t}%
^{v_{\cdot}},Z_{t}^{v_{\cdot}},k_{t}^{v_{\cdot}}  ,v_{t})  -\ell(
t,y_{t},Y_{t},z_{t},Z_{t},k_{t}  ,u_{t})   \big)
dt \big] .
\end{eqnarray*}
Since $\beta$\ and $\gamma$\ are convex, we obtain
\begin{eqnarray*}
&& \beta(  y_{T}^{v_{\cdot}})  -\beta(  y_{T}) \geq \langle\beta
_{y}(  y_{T}) ,   y_{T}^{v_{\cdot}}-y_{T}\rangle ,\\
&& \gamma(  Y_{0}^{v_{\cdot}})  -\gamma(  Y_{0})  \geq \langle\gamma
_{Y}(  Y_{0}),    Y_{0}^{v_{\cdot}}-Y_{0}\rangle,
\end{eqnarray*}
which imply that
\begin{eqnarray*}
&& \hspace{-1cm} J(  v_{\cdot})  -J(  u_{\cdot}) \geq\mathbb{E} \big[ \langle \beta
_{y}(  y_{T}),   y_{T}^{v_{\cdot}}-y_{T}\rangle \big]
+\mathbb{E} \big[ \langle \gamma_{Y}(  Y_{0}),    Y_{0}^{v_{\cdot}} -Y_{0}\rangle \big]  \\
&& + \; \mathbb{E} \big[ \int_{0}^{T} \big(   \ell(  t,y_{t}^{v_{\cdot}},Y_{t}^{v_{\cdot}},z_{t}^{v_{\cdot}},Z_{t}^{v_{\cdot}},k_{t}^{v_{\cdot}}  ,v_{t})  -\ell(
t,y_{t},Y_{t},z_{t},Z_{t},k_{t}  ,u_{t})   \big)
dt \big] .
\end{eqnarray*}
But from the adjoint equation (\ref{eq:3.1}) and system (\ref{eq:2.1}) we know
\begin{align*}
p_{0}=-\gamma_{Y}( Y_{0}) ,P_{T}=-c \, R^{\ast}  p_{T}+\beta_{y}( y_{T}) .
\end{align*}
Thus (\ref{eq:4.1}) holds.
\end{proof}

\bigskip

The following lemma contains duality relations between (\ref{eq:2.1}) and (\ref{eq:3.2}) (see the equivalent equations in the proof of
Theorem~\ref{thm:existence for adjoint eqn}).
\begin{lemma}\label{proof-lemma2}
Suppose that assumptions of Theorem~\ref{Theorem.3.1.} (in particular (A1)--(A3)) hold. Then
\begin{eqnarray}\label{eq:4.2}
&& -\; \mathbb{E} \big[ \left\langle  p_{0},  Y_{0}^{v_{\cdot}}-Y_{0}\right\rangle  \big]  =-\mathbb{E} \big[ \left\langle  p_{T},  Y_{T}^{v_{\cdot}}-Y_{T}\right\rangle  \big] \nonumber\\
&& - \; \mathbb{E} \big[ \int_{0}^{T} \left\langle   p_{t},f(  t,y_{t}^{v_{\cdot}},Y_{t}^{v_{\cdot}},z_{t}
^{v_{\cdot}},Z_{t}^{v_{\cdot}},k_{t}^{v_{\cdot}}  ,v_{t})  -f(
t,y_{t},Y_{t},z_{t},Z_{t},k_{t}  ,u_{t})
\right\rangle dt \big]\nonumber\\
&& + \; \mathbb{E} \, \big[ \int_{0}^{T}\left\langle H_{Y}(  t,y_{t},Y_{t},z_{t},Z_{t},k_{t},u_{t},p_{t},P_t,q_{t},Q_{t},V_{t}), Y_{t}^{v_{\cdot}}
-Y_{t}\right\rangle  dt \big] \nonumber\\
&& - \; \mathbb{E} \, \big[ \int_{0}^{T} \left\langle q_{t},g(  t,y_{t}^{v_{\cdot}},Y_{t}^{v_{\cdot}},z_{t}
^{v_{\cdot}},Z_{t}^{v_{\cdot}},k_{t}^{v_{\cdot}},v_{t}  ) -g(  t,y_{t}
,Y_{t},z_{t},Z_{t},k_{t},u_{t}  )
\right\rangle dt \big] \nonumber\\
&& + \; \mathbb{E} \, \big[ \int_{0}^{T}\left\langle H_{Z}(  t,y_{t},Y_{t},z_{t},Z_{t},k_{t} ,u_{t},p_{t},P_t,q_{t},Q_{t},V_{t}), Z_{t}^{v_{\cdot}}
-Z_{t}\right\rangle  dt \big] \nonumber\\
&& + \; \mathbb{E} \, \big[ \int_{0}^{T}\int_{\Theta}\langle H_{k}(  t,y_{t},Y_{t},z_{t}
,Z_{t},k_{t}  ,u_{t},p_{t},P_t,q_{t},Q_{t},V_{t}),\nonumber\\
&& \hspace{3in} k_{t}^{v_{\cdot}}(  \rho)  -k_{t}(  \rho)\rangle
\nu(  d\rho)  dt \big] ,
\end{eqnarray}
and
\begin{eqnarray}\label{eq:4.3}
&& \hspace{-1cm} \mathbb{E} \big[ \left\langle  P_{T},  y_{T}^{v_{\cdot}}-y_{T}\right\rangle  \big]  =\mathbb{E} \big[ \int_{0}^{T} \langle P_{t}, b(  t,y_{t}^{v_{\cdot}},Y_{t}^{v_{\cdot}},z_{t}^{v_{\cdot}}
,Z_{t}^{v_{\cdot}},k_{t}^{v_{\cdot}}  ,v_{t}) \nonumber \\
&& \hspace{2.3in} - \, b(
t,y_{t},Y_{t},z_{t},Z_{t},k_{t}  ,u_{t})
\rangle dt \big]\nonumber\\
&& + \, \mathbb{E}\big[ \int_{0}^{T}\left\langle H_{y}(  t,y_{t},Y_{t},z_{t},Z_{t},k_{t},u_{t},p_{t},P_t,q_{t},Q_{t},V_{t}),    y_{t}^{v_{\cdot}}%
-y_{t}\right\rangle  dt \big]\nonumber \\
&& + \, \mathbb{E} \big[ \int_{0}^{T}\left\langle H_{z}(  t,y_{t},Y_{t},z_{t},Z_{t},k_{t},u_{t},p_{t},P_t,q_{t},Q_{t},V_{t}),    z_{t}^{v_{\cdot}}
-z_{t}\right\rangle  dt \big]\nonumber\\
&& + \, \mathbb{E} \big[ \int_{0}^{T} \left\langle Q_{t}, \sigma(  t,y_{t}^{v_{\cdot}},Y_{t}^{v_{\cdot}},z_{t}^{v_{\cdot}},Z_{t}^{v_{\cdot}},k_{t}^{v_{\cdot}},v_{t}  )  -\sigma(
t,y_{t},Y_{t},z_{t},Z_{t},k_{t},u_{t}  )
\right\rangle dt \big]\nonumber\\
&& + \, \mathbb{E}\big[ \int_{0}^{T}\int_{\Theta} \langle V_{t}( \rho), \varphi(  t,y_t
^{v_{\cdot}},Y_t^{v_{\cdot}},z_{t}^{v_{\cdot}},Z_{t}^{v_{\cdot}},k_{t}^{v_{\cdot}},v_{t},\rho) \nonumber \\
 && \hspace{2in}  - \, \varphi(  t,y_t,Y_t,z_{t},Z_{t},k_{t},u_{t},\rho)  \rangle
\nu(  d\rho)  dt \big].
\end{eqnarray}
\end{lemma}
\begin{proof}
Applying integration by parts (Lemma~\ref{Lemma 4.1.}) to $\left\langle p_{t},Y_{t}^{v_{\cdot}}-Y_{t}\right\rangle $ gives
\begin{eqnarray*}
&& \hspace{-1cm} \left\langle  p_{T},  Y_{T}^{v_{\cdot}}-Y_{T}\right\rangle   = \left\langle  p_{0},  Y_{0}^{v_{\cdot}}-Y_{0}\right\rangle \\
&& - \; \int_{0}^{T} \left\langle  p_{t}, f(  t,y_{t}^{v_{\cdot}},Y_{t}^{v_{\cdot}},z_{t}
^{v_{\cdot}},Z_{t}^{v_{\cdot}},k_{t}^{v_{\cdot}}  ,v_{t})  -f(
t,y_{t},Y_{t},z_{t},Z_{t},k_{t}  ,u_{t})
\right\rangle dt\\
&& -\int_{0}^{T} \langle  p_{t}, \big( g(  t,y_{t}^{v_{\cdot}},Y_{t}^{v_{\cdot}},z_{t}
^{v_{\cdot}},Z_{t}^{v_{\cdot}},k_{t}^{v_{\cdot}}  ,v_{t})\\
&& \hspace{2.6in} -g(t,y_{t},Y_{t},z_{t},Z_{t},k_{t}  ,u_{t})\big) \overleftarrow{dB}_{t}
\rangle\\
&&+\int_{0}^{T} \int_{\Theta} \left\langle  p_{t}, (k_{t}^{v_{\cdot}}(\rho )-k_{t}(\rho))
\right\rangle \tilde{N}(  d\rho,dt)\\
&& +\int_{0}^{T} \left\langle  p_{t}, (Z_{t}^{v_{\cdot}}-Z_{t})dW_{t}
\right\rangle -\int_{0}^{T} \left\langle  Y_{t}^{v_{\cdot}}-Y_{t}, q_{t}\overleftarrow{dB}_{t}
\right\rangle\\
&& +\int_{0}^{T}\left\langle  H_{Y}(  t,y_{t},Y_{t},z_{t},Z_{t},k_{t},u_{t},p_{t},P_t,q_{t},Q_{t},V_{t}),Y_{t}^{v_{\cdot}}
-Y_{t}\right\rangle  dt\\
&& +\int_{0}^{T}\left\langle Y_{t}^{v_{\cdot}}
-Y_{t}, H_{Z}(  t,y_{t},Y_{t},z_{t},Z_{t},k_{t},u_{t},p_{t},P_t,q_{t},Q_{t},V_{t})dW_{t}\right\rangle \\
&& +\int_{0}^{T} \int_{\Theta} \langle Y_{t}^{v_{\cdot}}
-Y_{t}, H_{k}(  t,y_{t},Y_{t},z_{t},Z_{t},k_{t},u_{t},\\
&& \hspace{2.6in} p_{t},P_t,q_{t},Q_{t},V_{t})\rangle \tilde{N}(  d\rho,dt)\\
&& -\int_{0}^{T} \left\langle  q_{t},g(  t,y_{t}^{v_{\cdot}},Y_{t}^{v_{\cdot}},z_{t}
^{v_{\cdot}},Z_{t}^{v_{\cdot}},k_{t}^{v_{\cdot}},v_{t}  ) -g(  t,y_{t}
,Y_{t},z_{t},Z_{t},k_{t},u_{t}  )
\right\rangle dt\\
&& +\int_{0}^{T}\left\langle H_{Z}(  t,y_{t},Y_{t},z_{t},Z_{t},k_{t},u_{t},p_{t},P_t,q_{t},Q_{t},V_{t}), Z_{t}^{v_{\cdot}}
-Z_{t}\right\rangle  dt\\
&& +\int_{0}^{T}\int_{\Theta}\langle H_{k}(  t,y_{t},Y_{t},z_{t}
,Z_{t},k_{t}  ,u_{t},p_{t},P_t,q_{t},Q_{t},V_{t}),\\
&& \hspace{2.6in} k_{t}^{v_{\cdot}}(  \rho)  -k_{t}(  \rho)\rangle
\nu(  d\rho)  dt.
\end{eqnarray*}

Now by taking the expectation to the above equality, we obtain (\ref{eq:4.2}).

Similarly
\begin{eqnarray*}
&& \hspace{-1cm} \left\langle  P_{T},  y_{T}^{v_{\cdot}}-y_{T}\right\rangle    =\int_{0}^{T} \langle P_{t}, b(  t,y_{t}^{v_{\cdot}},Y_{t}^{v_{\cdot}},z_{t}^{v_{\cdot}}
,Z_{t}^{v_{\cdot}},k_{t}^{v_{\cdot}}  ,v_{t}) \\
&& \hspace{1.7in} - \; b(
t,y_{t},Y_{t},z_{t},Z_{t},k_{t}  ,u_{t})
\rangle dt\\
&& +\int_{0}^{T} \langle  P_{t}, (\sigma(  t,y_{t}^{v_{\cdot}},Y_{t}^{v_{\cdot}},z_{t}
^{v_{\cdot}},Z_{t}^{v_{\cdot}},k_{t}^{v_{\cdot}}  ,v_{t})\\
&& \hspace{1.7in} - \; \sigma(
t,y_{t},Y_{t},z_{t},Z_{t},k_{t}  ,u_{t}))dW_{t}
\rangle\\
&& +\int_{0}^{T} \int_{\Theta} \left\langle y_{t}^{v_{\cdot}}-y_{t},V_{t}(\rho)
\right\rangle \tilde{N}(  d\rho,dt)\\
&& -\int_{0}^{T} \left\langle  P_{t}, (z_{t}^{v_{\cdot}}-z_{t})\overleftarrow{dB}_{t}
\right\rangle +\int_{0}^{T} \left\langle y_{t}^{v_{\cdot}}-y_{t},Q_{t}dW_{t}
\right\rangle\\
&& +\int_{0}^{T}\int_{\Theta} \langle  P_{t}, \big(\varphi(  t,y_t
^{v_{\cdot}},Y_t^{v_{\cdot}},z_{t}^{v_{\cdot}},Z_{t}^{v_{\cdot}},k_{t}^{v_{\cdot}}  ,v_{t}
,\rho) \\
&& \hspace{1.7in}  - \; \varphi(  t,y_t,Y_t,z_{t},Z_{t},k_{t},u_{t},\rho)\big) \rangle  \tilde{N}(  d\rho,dt)
\\
&& +\int_{0}^{T}\left\langle  y_{t}^{v_{\cdot}}
-y_{t}, H_{z}(  t,y_{t},Y_{t},z_{t},Z_{t},k_{t},u_{t},p_{t},P_t,q_{t},Q_{t},V_{t})\overleftarrow{dB}_{t} \right\rangle \\
&& +\int_{0}^{T}\left\langle H_{y}(  t,y_{t},Y_{t},z_{t},Z_{t},k_{t},u_{t},p_{t},P_t,q_{t},Q_{t},V_{t}),    y_{t}^{v_{\cdot}}%
-y_{t}\right\rangle  dt \\
&& +\int_{0}^{T}\left\langle H_{z}(  t,y_{t},Y_{t},z_{t},Z_{t},k_{t},u_{t},p_{t},P_t,q_{t},Q_{t},V_{t}),    z_{t}^{v_{\cdot}}
-z_{t}\right\rangle  dt \\
&& +\int_{0}^{T} \langle Q_{t}, \sigma(  t,y_{t}^{v_{\cdot}},Y_{t}^{v_{\cdot}},z_{t}^{v_{\cdot}},Z_{t}^{v_{\cdot}},k_{t}^{v_{\cdot}},v_{t}  ) \\
&& \hspace{1.7in}-\sigma(
t,y_{t},Y_{t},z_{t},Z_{t},k_{t},u_{t}  )
\rangle dt \\
&& +\int_{0}^{T}\int_{\Theta} \langle V_{t}( \rho), \varphi(  t,y_t
^{v_{\cdot}},Y_t^{v_{\cdot}},z_{t}^{v_{\cdot}},Z_{t}^{v_{\cdot}},k_{t}^{v_{\cdot}}  ,v_{t}
,\rho) \\
&& \hspace{1.7in}  -\varphi(  t,y_t,Y_t,z_{t},Z_{t},k_{t} ,u_{t},\rho)  \rangle
\nu(  d\rho)  dt.
\end{eqnarray*}
By taking the expectation to this equality (\ref{eq:4.3}) holds.
\end{proof}

The remaining is devoted to completing the proof of Theorem~\ref{Theorem.3.1.}.

\bigskip

\noindent \begin{proof}[Proof of Theorem~\ref{Theorem.3.1.}]
Observe first from (\ref{eq:3.1}) that
\begin{eqnarray}\label{eq:4.4}
&& \hspace{-1cm}
\ell(  t,y_{t}^{v_{\cdot}},Y_{t}^{v_{\cdot}},z_{t}^{v_{\cdot}},
Z_{t}^{v_{\cdot}},k_{t}^{v_{\cdot}}  ,v_{t})  -\ell(t,y_{t},Y_{t},z_{t},Z_{t},k_{t}  ,u_{t}) \nonumber \\
&& = - \; \big(  H(  t,y_{t}^{v_{\cdot}},Y_{t}^{v_{\cdot}},z_{t}^{v_{\cdot}},Z_{t}^{v_{\cdot}}
,k_{t}^{v_{\cdot}}  ,v_{t},p_{t},q_{t},P_{t},Q_{t},V_{t}) \nonumber \\
&& \hspace{1.4in} - \; H(t,y_{t},Y_{t},z_{t},Z_{t},k_{t}  ,u_{t},p_{t},q_{t}
,P_{t},Q_{t},V_{t}) \big) \nonumber \\
&& \hspace{1cm} + \; \left\langle  p_{t}, f(  t,y_{t}^{v_{\cdot}},Y_{t}^{v_{\cdot}},z_{t}
^{v_{\cdot}},Z_{t}^{v_{\cdot}},k_{t}^{v_{\cdot}}  ,v_{t})  -f(
t,y_{t},Y_{t},z_{t},Z_{t},k_{t}  ,u_{t})
\right\rangle \nonumber \\
&& \hspace{1cm} - \;  \langle P_{t}, b(  t,y_{t}^{v_{\cdot}},Y_{t}^{v_{\cdot}},z_{t}^{v_{\cdot}}
,Z_{t}^{v_{\cdot}},k_{t}^{v_{\cdot}}  ,v_{t}) - b(t,y_{t},Y_{t},z_{t},Z_{t},k_{t}  ,u_{t}) \rangle \nonumber\\
&& \hspace{1cm} + \;  \langle q_{t}, g(  t,y_{t}^{v_{\cdot}},Y_{t}^{v_{\cdot}},z_{t}^{v_{\cdot}}
,Z_{t}^{v_{\cdot}},k_{t}^{v_{\cdot}}  ,v_{t}) - g(t,y_{t},Y_{t},z_{t},Z_{t},k_{t}  ,u_{t}) \rangle \nonumber\\
&& \hspace{1cm} - \; \left\langle Q_{t}, \sigma(  t,y_{t}^{v_{\cdot}},Y_{t}^{v_{\cdot}},z_{t}^{v_{\cdot}},Z_{t}^{v_{\cdot}},k_{t}^{v_{\cdot}},v_{t}  )
-\sigma(t,y_{t},Y_{t},z_{t},Z_{t},k_{t},u_{t}  )
\right\rangle \nonumber \\
&& \hspace{1cm} - \; \int_{\Theta} \langle V_{t}( \rho), \varphi(  t,y_t
^{v_{\cdot}},Y_t^{v_{\cdot}},z_{t}^{v_{\cdot}},Z_{t}^{v_{\cdot}},k_{t}^{v_{\cdot}},v_{t},\rho) \nonumber \\
&& \hspace{2in}  - \, \varphi(  t,y_t,Y_t,z_{t},Z_{t},k_{t},u_{t},\rho)  \rangle
\nu(  d\rho) .
\end{eqnarray}
Next apply Lemma~\ref{proof-lemma2} and (\ref{eq:4.4}) in Lemma~\ref{proof-lemma1} to find that
\begin{eqnarray}\label{eq:4.5}
&& \hspace{-0.75cm} J(  v_{\cdot})  - J(  u_{\cdot})  \nonumber \\
&& \hspace{-0.25cm} \geq \; \mathbb{E}\big[ \int_{0}^{T}\left\langle H_{y}(  t,y_{t},Y_{t},z_{t},Z_{t},k_{t} ,u_{t},p_{t},P_t,q_{t},Q_{t},V_{t}),    y_{t}^{v_{\cdot}}
-y_{t}\right\rangle  dt \big]\nonumber\\
&& \hspace{-0.25cm} + \; \mathbb{E}\big[ \int_{0}^{T}\left\langle H_{Y}(  t,y_{t},Y_{t},z_{t},Z_{t},k_{t},u_{t},p_{t},P_t,q_{t},Q_{t},V_{t}),   Y_{t}^{v_{\cdot}}
-Y_{t}\right\rangle  dt \big]\nonumber\\
&& \hspace{-0.25cm} + \; \mathbb{E}\big[ \int_{0}^{T}\left\langle H_{z}(  t,y_{t},Y_{t},z_{t},Z_{t},k_{t} ,u_{t},p_{t},P_t,q_{t},Q_{t},V_{t}),  z_{t}^{v_{\cdot}}
-z_{t}\right\rangle  dt \big]\nonumber\\
&& \hspace{-0.25cm} + \; \mathbb{E}\big[ \int_{0}^{T}\left\langle H_{Z}(  t,y_{t},Y_{t},z_{t},Z_{t},k_{t} ,u_{t},p_{t},P_t,q_{t},Q_{t},V_{t}),   Z_{t}^{v_{\cdot}}
-Z_{t}\right\rangle  dt \big]\nonumber\\
&& \hspace{-0.25cm} + \; \mathbb{E}\big[ \int_{0}^{T}\int_{\Theta}\langle H_{k}(  t,y_{t},Y_{t},z_{t}
,Z_{t},k_{t}  ,u_{t},p_{t},P_t,q_{t},Q_{t},V_{t}),  k_{t}^{v_{\cdot}}(  \rho) - \; k_{t}(  \rho)\rangle
\nu(  d\rho)  dt \big]\nonumber\\
&& \hspace{-0.25cm} - \; \mathbb{E}\big[ \int_{0}^{T}\big(  H(  t,y_{t}^{v_{\cdot}},Y_{t}^{v_{\cdot}},z_{t}^{v_{\cdot}},Z_{t}^{v_{\cdot}}
,k_{t}^{v_{\cdot}}  ,v_{t},p_{t},P_t,q_{t},Q_{t},V_{t}) \nonumber \\
&& \hspace{1cm} \hspace{1.4in} - \, H(t,y_{t},Y_{t},z_{t},Z_{t},k_{t}  ,u_{t},p_{t},P_t,q_{t},Q_{t},V_{t})\big)  dt \big].
\end{eqnarray}
Here we have used the formula $h(\omega , x) := c \,  R x + \xi (\omega) , x \in \mathbb{R}^n ,$ to get the cancelation
$$\mathbb{E} \big[ \left\langle   c \, R^{\ast} p_T , y_{T}^{v_{\cdot}}-y_{T}\right\rangle \big] -
\mathbb{E} \big[ \left\langle  p_T , Y_{T}^{v_{\cdot}}-Y_{T}\right\rangle  \big] = 0$$
resulting from (\ref{eq:4.1}) of Lemma~\ref{proof-lemma1} and (\ref{eq:4.3}) of Lemma~\ref{proof-lemma2}.

On the other hand, from the concavity condition (ii) of the mapping $$(y,Y,z,Z,k,v) \mapsto H(t,y,Y,z,Z,k,v,p_t,P_t,q_t,Q_t,V_t)$$ it follows that
\begin{eqnarray*}
&& \hspace{-0.5cm} H(  t,y_{t}^{v_{\cdot}},Y_{t}^{v_{\cdot}},z_{t}^{v_{\cdot}},Z_{t}^{v_{\cdot}},k_{t}^{v_{\cdot}}  ,v_{t}
,p_{t},P_t,q_{t},Q_{t},V_{t}) \nonumber \\
&& - \, H(  t,y_{t},Y_{t},z_{t},Z_{t}
,k_{t}  ,u_{t},p_{t},P_t,q_{t},Q_{t},V_{t}) \\
&& \leq \left\langle H_{y}(  t,y_{t},Y_{t},z_{t},Z_{t},k_{t}
,u_{t},p_{t},P_t,q_{t},Q_{t},V_{t}), y_{t}^{v_{\cdot}}-y_{t}\right\rangle \\
&& + \, \left\langle H_{Y}(  t,y_{t},Y_{t},z_{t},Z_{t},k_{t}
,u_{t},p_{t},P_t,q_{t},Q_{t},V_{t}),   Y_{t}^{v_{\cdot}}-Y_{t}\right\rangle \\
&& + \, \left\langle H_{z}(  t,y_{t},Y_{t},z_{t},Z_{t},k_{t}
,u_{t},p_{t},P_t,q_{t},Q_{t},V_{t}),   z_{t}^{v_{\cdot}}-z_{t}\right\rangle \\
&& + \, \left\langle H_{Z}(  t,y_{t},Y_{t},z_{t},Z_{t},k_{t}
,u_{t},p_{t},P_t,q_{t},Q_{t},V_{t}),   Z_{t}^{v_{\cdot}}-Z_{t}\right\rangle \\
&& +\int_{\Theta}\left\langle H_{k}(  t,y_{t},Y_{t},z_{t},Z_{t},k_{t},u_{t},p_{t},P_t,q_{t},Q_{t},V_{t}),   k_{t}^{v_{\cdot}}(\rho)  -k_{t}(  \rho)\right\rangle    \nu(  d\rho) \\
&& + \, \left\langle H_{v}(  t,y_{t},Y_{t},z_{t},Z_{t},k_{t}
,u_{t},p_{t},P_t,q_{t},Q_{t},V_{t}),  v_{t}-u_{t}\right\rangle  .
\end{eqnarray*}
In particular,
\begin{eqnarray*}
&& \hspace{-1cm} - \, \left\langle H_{v}(  t,y_{t},Y_{t},z_{t},Z_{t},k_{t}
,u_{t},p_{t},P_t,q_{t},Q_{t},V_{t}),  v_{t}-u_{t}\right\rangle \\
&& \leq \left\langle H_{y}(  t,y_{t},Y_{t},z_{t},Z_{t},k_{t}
,u_{t},p_{t},P_t,q_{t},Q_{t},V_{t}),  y_{t}^{v_{\cdot}}-y_{t}\right\rangle \\
&& + \, \left\langle H_{Y}(  t,y_{t},Y_{t},z_{t},Z_{t},k_{t}
,u_{t},p_{t},P_t,q_{t},Q_{t},V_{t}),  Y_{t}^{v_{\cdot}}-Y_{t}\right\rangle \\
&& + \, \left\langle H_{z}(  t,y_{t},Y_{t},z_{t},Z_{t},k_{t}
,u_{t},p_{t},P_t,q_{t},Q_{t},V_{t}),  z_{t}^{v_{\cdot}}-z_{t}\right\rangle \\
&& + \, \left\langle H_{Z}(  t,y_{t},Y_{t},z_{t},Z_{t},k_{t}
,u_{t},p_{t},P_t,q_{t},Q_{t},V_{t}),  Z_{t}^{v_{\cdot}}-Z_{t}\right\rangle \\
&& + \, \int_{\Theta}\left\langle H_{k}(  t,y_{t},Y_{t},z_{t},Z_{t},k_{t} ,u_{t},p_{t},P_t,q_{t},Q_{t},V_{t}),  k_{t}^{v_{\cdot}}(
\rho)  -k_{t}(  \rho)  \right\rangle  \nu(  d\rho) \\
&& - \, [   H(  t,y_{t}^{v_{\cdot}},Y_{t}^{v_{\cdot}},z_{t}^{v_{\cdot}},Z_{t}^{v_{\cdot}},k_{t}^{v_{\cdot}}
,v_{t},p_{t},P_t,q_{t},Q_{t},V_{t}) \nonumber \\ && \hspace{1.75in} - \, H(  t,y_{t},Y_{t},z_{t}
,Z_{t},k_{t}  ,u_{t},p_{t},P_t,q_{t},Q_{t},V_{t})
 ] .
\end{eqnarray*}

Now by applying this latter result in (\ref{eq:4.5}) we obtain
\begin{equation}\label{eq:4.6}
J(  v_{\cdot})  -J(  u_{\cdot})  \geq- \mathbb{E}\big[ \int_{0}^{T}\left\langle H_{v}(
t,y_{t},Y_{t},z_{t},Z_{t},k_{t}  ,u_{t},p_{t},P_t,q_{t},Q_{t},V_{t}),   v_{t}-u_{t}\right\rangle  dt \big].
\end{equation}

On the other hand, the maximum condition (iii) yields
\begin{equation*}
\left\langle H_{v}( t,y_{t},Y_{t},z_{t},Z_{t},k_{t},u_{t},p_{t},P_t,q_{t},Q_{t},V_{t}),  v_{t}-u_{t}\right\rangle \leq 0 .
\end{equation*}
Hence (\ref{eq:4.6}) becomes
\[
J( v_{\cdot}) -J( u_{\cdot}) \geq 0 .
\]
Since $u_{\cdot}$ is an arbitrary element of $\mathcal{U}_{ad},$ this inequality completes the proof if we recall (\ref{eq:2.3}).
\end{proof}

\bigskip

{\bf Acknowledgement.}
The authors would like to thank the associate editor and anonymous referee(s) for their remarks, which have helped in improving the first version of this paper.

\fussy


\begin{thebibliography}{99}
\bibitem {Ant93} F. Antonelli, Forward-backward stochastic differential
equations, Ann. Appl. Probab. 3 (1993), 777--793.
\bibitem {Bah-Mez09} K. Bahlali, N. Khelfallah and B. Mezerdi,
Necessary and sufficient conditions for near-optimality in stochastic control of FBSDEs,
Systems Control Lett. 58 (2009), 12, 857--864.
\bibitem {Ming-Max09} Q. Meng, Optimal control problem of fully coupled
forward-backward stochastic systems with Poisson jumps under partial
information, arXiv:0911.3225v1 [math.OC], 2009.
\bibitem{Nua-Pard} D. Nualart and E. Pardoux, Stochastic calculus associated with Skokrohod's integral,
Lecture Notes in Control and Inf. Sciences, 96 (1987), 363--372.
\bibitem {Oks-Sulem08} B. {\O}ksendal and A. Sulem, Maximum principles for
optimal control of forward-backward stochastic differential equations with
jumps, SIAM J. Control Optim., 48 (2009), 5, 2945--2976.
\bibitem {PP2} E. Pardoux and S. Peng, Backward doubly stochastic differential
equations and system of quasilinear SPDEs,  Prob. Th. \& Rel Fields, 98
(1994), no. 2, 209--227.
\bibitem {PS} S. Peng and Y. Shi, A type-symmetric forward-backward stochastic
differential equations, C. R. Acad. Sci. Paris Ser. I, 336 (2003), no. 1, 773--778.
\bibitem {Pe-Wu99} S. Peng and Z. Wu, Fully coupled forward-backward stochastic
differential equations and applications to optimal control, SIAM J. Control
Optim, 37 (1999), 825--843.
\bibitem {ZS2} Zhu Qingfeng and Yufeng Shi, Forward-backward doubly stochastic differential equations with
random jumps and stochastic partial differential-integral equations, [OL].[201001-1044],
\url{http://www.paper.edu.cn/index.php/default/en-releasepaper/downPaper/201001-1044.}
\bibitem {Shi2012} J. T. Shi, Necessary conditions for optimal control of
forward-backward stochastic systems with random jumps, International Journal
of Stochastic Analysis Volume 2012, Article ID 258674, doi:10.1155/2012/258674.
\bibitem {Shi-Wu06} J. T. Shi and Z. Wu, The maximum principle for fully coupled
forward-backward stochastic control system, Acta Automatica Sinica, 32 (2006), 2, 161--169.
\bibitem{Jiang-Xu} Xu Shuli and Jiang Jun,
Maximum principles for forward-backward doubly stochastic differential equations with jumps, [OL],
\url{http://www.paper.edu.cn/en_releasepaper/content/4486654.}
\bibitem {Situ} R. Situ, On solution of backward stochastic differential equations with jumps and
applications, Stoch. Process. Appl., 66 (1997), 209-–236.
\bibitem {SL} X. Sun and Y. Lu, The property for solutions of the
multi-dimensional backward doubly stochastic differential equations with
jumps, Chin. J. Appl. Probab. Stat., 24 (2008), 73--82.
\bibitem {WW} X. Wang and Z. Wu, FBSDE with Poisson process and its application
to linear quadratic stochastic optimal control problem with random jumps, Acta
Automatica Sinica, 29 (2003), 821--826.
\bibitem {Wu98} Z. Wu, Maximum principle for optimal control problem of fully
coupled forward-backward stochastic systems, Systems Sci. Math. Sci., 11 (1998), 3, 249--259.
\bibitem {W2} Z. Wu, Forward-backward stochastic differential equations with
Brownian motion and Poisson process, Acta Math. Appl. Sinica, 15 (1999), 433--443.
\bibitem {Xu95} W. Xu, Stochastic maximum principle for optimal control problem
of forward and backward system, J. Australian Mathematical Society B, 37
(1995), 172--185.
\bibitem {YS} Juliang Yin and R. Situ, On solutions of forward-backward stochastic
differential equations with Poisson jumps, Stoch. Anal. Appl., 21 (2003), 1419--1448.
\bibitem{Yong97} J. Yong, Finding adapted solutions of forward-backward stochastic
differential equations--method of continuation, Prob. Th. \& Rel Fields, 107 (1997), 537--572.
\bibitem {Zha-Shi2012} Liangquan Zhang and Yufeng Shi, General doubly stochastic maximum
principle and its applications to optimal control of stochastic partial differential
equations, arXiv:1009.6061v3[math.OC], 2012.
\end{thebibliography}
\end{document}